\documentclass[10pt]{article}
\def\date{4.4.12} %
\usepackage{hyperref}
\usepackage{amssymb, amsmath}
\input{liemacs10.sty} 

\newcommand\uline{\underline}
 
\newcommand\cR{\mathcal{R}} 
\newcommand\co{\mathop{\rm co}\nolimits}

\renewcommand{\mlabel}{\label}

\begin{document} 

\title{On convex hulls of orbits of Coxeter groups\\ and Weyl groups} 
\author{G. Hofmann\begin{footnote}{
Department of Mathematics and Statistics, Chase Building, Dalhousie
University, Halifax, Nova Scotia  B3H 3J5,  Canada, 
hofmann@mathstst.dal.ca}\end{footnote} and 
K.-H. Neeb\begin{footnote}{
Department  Mathematik, FAU Erlangen-N\"urnberg, Cauerstrasse 11, 
91058 Erlangen, Germany, karl-hermann.neeb@math.uni-erlangen.de}
\end{footnote}
\begin{footnote}{Supported by DFG-grant NE 413/7-1, Schwerpunktprogramm 
``Darstellungstheorie''.} 
\end{footnote}}

\maketitle

\begin{abstract} The notion of a linear Coxeter system introduced by Vinberg generalizes the geometric 
representation of a Coxeter group. Our main theorem asserts that if $v$ is an element of the 
Tits cone of a linear Coxeter system and $\cW$ is the corresponding Coxeter group, then 
$\cW v \subeq v - C_v,$ where $C_v$ is the convex cone generated by 
the coroots $\check \alpha$, for which $\alpha(v) > 0$. This implies that the 
convex hull of $\cW v$ is completely determined by the image of $v$ under the reflections 
in $\cW$. We also apply an analogous result for convex hulls of $\cW$-orbits in the dual 
space, although this action need not correspond to a linear Coxeter system. 
Motivated by the applications in representation theory, we further 
extend these results to Weyl group orbits of locally finite 
and locally affine root systems. In the locally affine case, we also derive 
some applications on minimizing linear functionals on Weyl group orbits. 
\end{abstract}

\section*{Introduction} 

The present paper is motivated by the unitary representation theory 
of locally finite, resp., locally affine Lie algebras and their 
analytic counterparts (\cite{Ne98, Ne10, Ne12}). 
In the algebraic context, these Lie algebras  
$\g$ contain a maximal abelian subalgebra $\ft$ for which the complexification 
$\g_\C$ has a root decomposition 
$\g_\C = \ft_\C \oplus \bigoplus_{\alpha \in \Delta} \g_\C^\alpha$   
and there is a distinguished class of unitary representations 
$(\rho,V)$ of $\g$ on a pre-Hilbert space $V$ for which 
the operators $\rho(x)$, $x \in \g$, are skew-symmetric, 
on which $\rho(\ft)$ is diagonalizable,  
and the corresponding weight set $\cP_V \subeq i \ft^*$ has the form 
\[ \cP_V = \conv(\cW \lambda) \cap (\lambda + \cQ),\] 
where $\cW \subeq \GL(\ft)$ is the {\it Weyl group} of the pair 
$(\g,\ft)$, and $\cQ \subeq i\ft^*$ is the {\it root group}. Then 
\[ \Ext(\conv(\cP_V)) = \cW \lambda \] 
is the set of {\it extremal weights}. For finite-dimensional 
compact Lie algebras $\fg$ and the unitary (=compact) forms of Kac--Moody algebras, 
this is the well-known description of the weight set of unitary highest weight 
modules (\cite{Ka90}).  For the generalization to the locally finite, resp., locally affine 
case we refer to \cite{Ne98}, resp., \cite{Ne10} for details. 
Motivated by these representation theoretic issues, 
the present paper addresses a better understanding 
of the convex hulls of Weyl group orbits in $\ft$ and $\ft^*$. 

In the locally finite, resp., locally affine case, 
the Weyl group is a direct limit of finite, resp., affine Coxeter groups, 
and it is the geometry of linear actions of Coxeter groups that provides the 
key to the crucial information on convex hulls of orbits. 
Therefore this paper is divided into two part. 
The first part deals with {\it linear Coxeter systems} 
(\cite{Vin71}). These are realizations of Coxeter groups by groups 
generated by a finite set 
$(r_s)_{s \in S}$ of reflections of a finite-dimensional vector space~$V$ 
satisfying the following conditions. 
We write $r_s(v) = v - \alpha_s(v) \check \alpha_s$ with 
elements $\alpha_s \in V^*$ and $\check \alpha_s \in V$ and 
$\cone(M)$ for the convex cone generated by the subset $M$ of a real vector space. 
Then the conditions for a linear Coxeter system are: 
\begin{description}
\item[\rm(LCS1)] The polyhedral cone 
$K := \{ v \in V \: (\forall s \in S)\, \alpha_s(v) \geq 0 \}$ 
has interior points. 
\item[\rm(LCS2)] $(\forall s \in S)\, \alpha_s \not\in \cone(\{ \alpha_t \: 
t\not=s\})$. 
\item[\rm(LCS3)] $(\forall w \in \cW\setminus \{\1\})
\, wK^0 \cap K^0 =\eset$, where $K^0$ denotes the interior of $K$. 
\end{description}

For any linear Coxeter system, the set $T := \cW K\subeq V$ is a convex cone,  
called the {\it Tits cone}. A {\it reflection in $\cW$} is an element 
conjugate to one of the $r_s$, $s \in S$. Any reflection 
can be written as $r_\alpha(v) = v - \alpha(v) \check \alpha$ 
with $\alpha \in V^*$ and $\check \alpha \in V$, where 
$\alpha$ belongs to the set 
$\Delta := \cW \{ \alpha_s \: s \in S \}$ of {\it roots} and 
$\check\alpha \in \check\Delta := \cW \{ \check\alpha_s \: s \in S \}$ is a 
{\it coroot}. In these terms, our first main result 
(Theorem~\ref{thm:cox-cont}) asserts that 
\begin{equation}
  \label{eq:con1}
\cW v \subeq v - C_v \quad \mbox{ for } \quad 
C_v := \cone \{ \check \alpha \: \alpha(v) >  0 \} \quad \mbox{ and } \quad 
v \in T = \cW K.
\end{equation}
As an inspection of two-dimensional examples shows, 
the restriction to elements of the Tits cone is crucial and that 
the boundary $\partial T$ may contain elements $v$ for which 
$\cW v$ is not contained in $v - C_v$.\begin{footnote}
{The proof of Theorem~\ref{thm:cox-cont} relies on various 
results from the unpublished diploma thesis of Georg Hofmann \cite{HoG99} which contains already a proof for the 
special case where $v \in K^0$, i.e., where the stabilizer of $v$ is trivial.}
\end{footnote}

For the applications to orbits of weights, we also 
need a corresponding result for elements in the dual space. 
Here a difficulty arises from the fact that the concept of a linear 
Coxeter system is not preserved by exchanging the role of $V$ and $V^*$, 
so that Theorem~\ref{thm:cox-cont} cannot be applied directly. However, this can be overcome 
by reduction to the subspace $U \subeq V^*$ 
generated by the dual cone $C_S^\star$ of $C_S := \cone\{ \check \alpha_s \: s \in S \}$. 
Then $\cW C_S^\star$ is the Tits cone for a linear Coxeter system on $U$ 
(Theorem~\ref{thm:cox-cont2}). 

{\bf Structure of this paper:} In Section~\ref{sec:1} we recall 
some basics on linear Coxeter systems and generalize some results 
which are well-known for the geometric representation of a Coxeter group 
to linear Coxeter systems. Here the main results are 
Theorem~\ref{thm:d.7} relating positivity conditions to the length 
function and Proposition~\ref{prop:d.8} asserting that the stabilizer 
of an element in the fundamental chamber $K$ is generated by reflections. 
Our main results,  \eqref{eq:con1} and its dual version, 
are proved in Section~\ref{sec:2}. 
In Section~\ref{sec:3} we turn to the applications to Weyl group orbits 
of linear functionals for locally finite and locally affine root systems 
(Theorems~\ref{thm:cox-cont} and \ref{thm:cox-cont2}). 
The two final subsections are also motivated by the representation 
theoretic problem to determine maximal and minimal eigenvalues of 
elements of the Cartan subalgebra in extremal weight representations. 
In this context we provide in the locally affine case 
a complete classification of the set 
$\cP_d^+$ of those linear functionals 
$\lambda$ for which $\lambda(d) = \min (\hat\cW\lambda)(d)$ holds  
for a certain distinguished element $d$ and show that 
any orbit of the affine Weyl group $\hat\cW$ intersects $\cP_d^+$ in an orbit of the 
corresponding locally finite Weyl group $\cW$. 
Writing $\hat\cW$ as a semidirect product 
$\cN \rtimes \cW$, where $\cN$ acts by unipotent isometries of a 
Lorentzian form, the minimization of the $d$-value on a 
$\hat\cW$-orbit turns into a problem of minimizing a quadratic form 
on an infinite-dimensional analog of a lattice in a euclidean space. 
These issues are briefly discussed in Subsection~\ref{subsec:3.4}. 

The results of the present paper constitute a crucial ingredient in the 
classification of semibounded unitary representations 
of double extensions of loop groups with values in Hilbert--Lie groups 
carried out in \cite{Ne12}.

{\bf Notation:} For a subset $E$ of the real vector space $V$ we 
write $\conv(E)$ for its convex hull and 
$\cone(E)$ for the convex cone generated by $E$. 
For a convex cone $C \subeq V$ we write 
$H(C) := C \cap - C$ for the largest subspace contained in $C$. 
We also write 
\[ E^\star := \{ \alpha \in V^* \: (\forall v \in E)\, \alpha(v) \geq 0 \} \] 
for its dual cone in $V^*$. For $E \subeq V^*$ we  define its dual 
cone by 
\[ E^\star := \{ v \in V \: (\forall \alpha \in E) \, \alpha(v) \geq 0 \}. \]

\tableofcontents 

\section{Linear Coxeter systems} 
\mlabel{sec:1}

In this section we recall Vinberg's concept of a linear Coxeter system, generalizing the  
geometric representation of a Coxeter group. 

\begin{defn} (a) 
Let $V$ be a real vector space. 
A {\it reflection data on $V$} consists of a family 
$(\alpha_s)_{s \in S}$ of linear functionals on $V$ 
and a family $(\check\alpha_s)_{s \in S}$ of elements of $V$ 
satisfying 
\[\alpha_s(\check \alpha_s) = 2 \quad \mbox{ for } \quad s \in S.\] 
Then $r_s(v) := v - \alpha_s(v) \check \alpha_s$ 
is a reflection on $V$. We write $\cW := \la r_s \: s \in S \ra 
\subeq \GL(V)$ for the subgroup generated by these reflections 
and 
\begin{equation}
  \label{eq:co}
 \co(v) :=  \conv(\cW v)
\end{equation}
for the convex hull of a $\cW$-orbit. 
We say that a reflection data is {\it of finite type} 
if $S$ is finite and $\dim V < \infty$. 

(b) We consider the following polyhedral cones 
in $V$ resp.~$V^*$: 
\[ C_S := \cone\{ \check \alpha_s \: s \in S \} \subeq V, \qquad 
\check C_S := \cone\{ \alpha_s \: s \in S \} \subeq V^*, \] 
and the {\it fundamental chamber} 
\[ K := \{ v \in V \: (\forall s \in S)\, \alpha_s(v) \geq 0 \}  
= (\check C_S)^\star.\]

(c) A reflection data of finite type is called a  {\it linear Coxeter system} 
(cf.\ \cite{Vin71}) if 
\begin{description}
\item[\rm(LCS1)] $K$ has interior points, i.e., 
the cone $\check C_S \subeq V^*$ is pointed. 
\item[\rm(LCS2)] $(\forall s \in S)\, \alpha_s \not\in \cone(\{ \alpha_t \: 
t\not=s\})$. 
\item[\rm(LCS3)] $(\forall w \in \cW\setminus \{\1\})
\, wK^0 \cap K^0 =\eset$. 
\end{description}
Then $T := \cW K$ is called the associated {\it Tits cone}. 
\end{defn}

\begin{defn} Let $V$ be a real vector space endowed with a symmetric 
bilinear form $\beta$. 
A {\it symmetric reflection data on $V$} consists of a family 
$(\check\alpha_s)_{s \in S}$ of non-isotropic elements of $V$ 
for which the linear functionals 
\[ \alpha_s(v) := \frac{2\beta(v, \check \alpha_s)}{\beta(\check \alpha_s, \check \alpha_s)} 
 \] 
define a reflection data on $V$. 
Then the reflections 
\[ r_s(v) 
= v - \alpha_s(v) \check \alpha_s 
= v - 2 \frac{\beta(v, \check \alpha_s)}{\beta(\check \alpha_s, \check \alpha_s)} \check \alpha_s  \] 
preserve $\beta$, so that $\cW \subeq \OO(V,\beta)$. 
\end{defn}

\begin{rem} \mlabel{rem:a.2} (a) Condition (LCS2) means that, for each $s \in S$, 
the dual cone of \break $\cone\{ \alpha_t \: t \not= s\}$ is strictly 
larger than $K$, i.e., there exists an element $v \in K\cap \ker \alpha_s$ 
with $\alpha_t(v) > 0$ for $t \not= s$. Then 
$K \cap \ker \alpha_s$ is a codimension $1$ face of the polyhedral cone~$K$. 

(b) Typical examples of linear Coxeter systems 
arise from the geometric representation of a 
Coxeter system $(\cW,S)$ (cf.\ \cite[\S 5.3]{Hu92}, \cite{Vin71}). 
\end{rem}

The following criterion 
for the recognition of linear Coxeter systems will be convenient 
in many situations because conditions (C1) and (C2) are preserved by 
the passage to the dual reflection data obtained by exchanging $V$ and~$V^*$. 

\begin{prop} \mlabel{prop:d.3} 
Let $(V, (\alpha_s)_{s \in S}, (\check \alpha_s)_{s \in S})$ 
be a reflection data of finite type satisfying (LCS1). 
Then it defines a linear Coxeter system if and only if the 
following conditions are satisfied for $s\not=t \in S$: 
\begin{description}
\item[\rm(C1)] $\alpha_s(\check \alpha_t)$ and 
$\alpha_t(\check \alpha_s)$ are either both negative or both 
zero. 
\item[\rm(C2)] $\alpha_s(\check \alpha_t)\alpha_t(\check \alpha_s) 
\geq 4$ or $\alpha_s(\check \alpha_t)\alpha_t(\check \alpha_s) 
= 4 \cos^2\frac{\pi}{k}$ for some natural number $k \geq 2$. 
\end{description}

In this case $(\cW, (r_s)_{s \in S})$ is a Coxeter system. 
\end{prop}

\begin{prf} 
According to 
\cite[p.~1085]{Vin71}, for a reflection data of finite type 
satisfying (LCS1/2), condition 
(LCS3) is equivalent to both conditions (C1/2). 

Therefore it remains to show that (C1) implies (LCS2). 
So we assume that 
$\alpha_s = \sum_{t \not=s }\lambda_t \alpha_t$ with finitely many 
non-zero $\lambda_t \geq 0$. Then (C1) leads to the contradiction 
\[ 2 = \alpha_s(\check \alpha_s) 
= \sum_{t \not=s }\lambda_t \alpha_t(\check \alpha_s) \leq 0. \] 
Therefore (C1) implies (LCS2). 
That $(\cW, (r_s)_{s \in S})$ is a Coxeter system follows from 
\cite[Thm.~2(6)]{Vin71}. 
\end{prf}

\begin{rem} \mlabel{rem:1.6} As a consequence of the preceding theorem, 
we obtain for every subset $S_0 \subeq S$ and every 
linear Coxeter system $(V, (\alpha_s)_{s \in S}, 
(\check \alpha_s)_{s \in S})$ a linear Coxeter system 
\[ (V, (\alpha_s)_{s \in {S_0}}, 
(\check \alpha_s)_{s \in {S_0}}).\] 
\end{rem}

\begin{rem}  \mlabel{rem:1.7} 
(a) If $(V, (\alpha_s)_{s \in S}, (\check \alpha_s)_{s \in S})$ is a reflection 
data, then we also consider the elements $\check \alpha_s$ as 
linear functionals on $V^*$. We thus obtain a reflection data 
$(V^*, (\check\alpha_s)_{s \in S}, (\alpha_s)_{s \in S})$. 
Suppose that $(V, (\alpha_s)_{s \in S}, 
(\check \alpha_s)_{s \in S})$ is a linear Coxeter system. 
Then (C1/2) also hold for the dual reflection data 
$(V^*, (\check\alpha_s)_{s \in S}, (\alpha_s)_{s \in S})$. 
Hence it is a linear Coxeter system if and only if 
(LCS1) holds, i.e., if the convex cone $C_S$ 
is pointed, i.e., $H(C_S) = \{0\}$. 

(b) If $C_S$ is not pointed, then we still obtain a linear Coxeter system 
by replacing $V^*$ by the smaller subspace 
\[ U := \Spann(C_S^\star) = H(C_S)^\bot \subeq V^*.\] 
Let $q \: V \to V/U^\bot \cong U^*$ 
denote the canonical projection. We put 
\[ \tilde S := \{ s \in S \: \check \alpha_s \not\in U^\bot = H(C_S)\}.\] 
Proposition~\ref{prop:d.3} implies that 
$(U, (q(\check \alpha_s))_{s \in \tilde S}, (\alpha_s)_{s \in \tilde S})$ 
is a linear Coxeter system because the convex cone 
$\check C_{\tilde S} \subeq \check C_S$ is pointed. 
We write 
\[ \cW_U  := \la r_s \: s \in \tilde S \ra \subeq \GL(U) \] 
for the corresponding reflection group. 

To relate this group to $\cW$, we first claim that $U$ is $\cW$-invariant. 
For $s \in S \setminus \tilde S$ the relation $\check\alpha_s \in U^\bot$ 
implies that the reflection $r_s$ acts trivially on the subspace 
$U \subeq V^*$. 
Next we observe that, if an element 
$\sum_{s \in S} \lambda_s \check \alpha_s \in C_S$ with 
$\lambda_s \geq 0$ is contained in $H(C_S)$, 
then $\lambda_s > 0$ implies that $\check \alpha_s \in H(C_S)$, 
i.e., $s \in S \setminus \tilde S$. We conclude that 
\[ H(C_S) = \cone \{ \check \alpha_s \: s \in S \setminus 
\tilde S \} = C_{S \setminus \tilde S}.\] 
For $s \in \tilde S$ we now derive from (C1) that 
$\alpha_s \in - (C_{S\setminus \tilde S})^\star = -H(C_S)^\star 
= U.$ 
Therefore $U$ is also invariant under $r_s$. 
Hence $U$ is $\cW$-invariant. As the reflections 
$r_s$, $s \in S \setminus \tilde S$, act trivially on $U$, 
we see that the restriction map 
\[ R\: \cW \to \cW_U, \quad w\mapsto w\res_U \] 
is a surjective homomorphism. 
\end{rem}

\begin{lem} \mlabel{lem:a.1} 
Let $a \geq 2$ and 
\[ r_1 := \pmat{-1 & a \\ 0 & 1}, \qquad 
   r_2 := \pmat{1 & 0 \\ a & -1} \in \GL_2(\R).\] 
If $e_1$ and $e_2$ denote the standard basis of $\R^2$, then 
\[ (r_1 r_2)^n e_1, r_2(r_1 r_2)^n e_1 \in \cone\{e_1, e_2\}
\quad \mbox{ for } \quad n \in \N_0.\] 
\end{lem}

\begin{prf} We define a linear functional on $\R^2$ by 
$\beta(x) := x_1 - \frac{1}{2}a x_2$. First we show by induction on $n$ 
that 
\[ (r_1 r_2)^n e_1 \in \cone\{e_1, e_2\} \quad \mbox{ and } \quad 
 \beta((r_1 r_2)^n e_1) \geq 0 \quad \mbox{ hold for any } \quad n \in \N_0.\] 
For $n = 0$ this is trivial. So 
let $v := (r_1 r_2)^n e_1$ and assume that the assertion holds for 
some $n \in \N_0$. With 
\[ r_1 r_2 = \pmat{ a^2 - 1 & - a\\ a & -1} \] 
we then obtain 
\begin{align*}
\beta(r_1 r_2 v) 
&= (a^2 - 1)v_1 - a v_2 - \frac{1}{2}a (av_1 - v_2) 
= (\frac{a^2}{2} - 1) v_1 - \frac{a}{2} v_2\\
&= \underbrace{(\frac{a^2}{2} - 2)}_{\geq 0} v_1 + \beta(v) \geq 0.
\end{align*}
We also have 
\[ (a^2 -1)v_1 - av_2 
= (a^2 - 3) v_1 + 2\Big(v_1 - \frac{a}{2}v_2\Big) 
\geq v_1 + 2\beta(v) \geq 0 \] and 
\begin{equation}
  \label{eq:e1}
av_1 - v_2 
= \frac{1}{a}\big((a^2 - 2)v_1 + 2(v_1 - \frac{a}{2}v_2)\big) 
\geq  \frac{1}{a}(2 v_1 + 2\beta(v))\geq 0.
\end{equation}
This proves that $r_1 r_2 v \in \cone\{e_1, e_2\}$ and completes our induction. 

For $v' := r_2 v = r_2 (r_1 r_2)^n e_1$, we finally obtain with \eqref{eq:e1} 
\[ v_1' = v_1 \geq 0 \quad \mbox{ and } \quad 
v_2' = av_1 - v_2 \geq 0,\] 
hence that $v' \in \cone\{e_1, e_2\}$.
\end{prf}

\begin{lem}
  \mlabel{lem:x} 
Let 
$(V, (\alpha_s, \alpha_t), (\check \alpha_s, \check \alpha_t))$ be a 
reflection data satisfying (C1/2). Then the subgroup 
$\Gamma \subeq \GL(V)$ generated by the reflections 
$r_s$ and $r_t$ is a Coxeter group. Its length function $\ell$ 
satisfies for $g \in \cW$: 
\[ \ell(gr_s) \geq \ell(g) \quad \Rarrow \quad g\check \alpha_s \in C_S 
= \cone \{ \check \alpha_s, \check \alpha_t \}.\] 
\end{lem}

\begin{prf} In view of $\ell(gr_s) \geq \ell(g)$, 
every reduced expression for $g$ is an alternating product of 
$r_s$ and $r_t$ ending in $r_t$. 
After normalizing the pair $(\alpha_s, \check \alpha_s)$ suitably, we 
may assume that $\alpha_s(\check \alpha_t) = \alpha_t(\check \alpha_s)$. 
Let 
\[ a := - \alpha_s(\check \alpha_t) \geq 0 \quad \mbox{ and } \quad 
U := \Spann \{ \check \alpha_s, \check \alpha_t\}.\] 
Then 
\[ \det\pmat{ 
\alpha_s(\check \alpha_s) & \alpha_s(\check \alpha_t) \\
\alpha_t(\check \alpha_s) & \alpha_t(\check \alpha_t)} 
= 4 - a^2 \] 
shows that, if $a^2 \not=4$, then 
$\check \alpha_s$ and $\check \alpha_t$ are linearly independent 
and the same holds for the restrictions of 
$\alpha_s$ and $\alpha_t$ to $U$. 

If $\check\alpha_s$ and $\check \alpha_t$ are linearly independent, then the 
corresponding matrices of the restrictions of $r_s$ and $r_t$ 
to $U$ are given by 
\[ r_1 := \pmat{-1 & a \\ 0 & 1} \quad \mbox{ and } \quad 
   r_2 := \pmat{1 & 0 \\ a & -1}.\] 
Moreover, $r_1$ and $r_2$ are orthogonal with respect to the 
symmetric bilinear form $\la \cdot, \cdot \ra$ defined with respect to the basis 
$(\check \alpha_s,\check \alpha_t)$ by the matrix 
\[ \pmat{ 1 & -\frac{a}{2} \\ 
-\frac{a}{2} & 1} \] 
because 
$\alpha_s(v) = 2 \la \check \alpha_s, v\ra$ and 
$\alpha_t(v) = 2 \la \check \alpha_t, v\ra$ hold for each $v \in U$. 

{\bf Case 1:} $a^2 < 4$. In this case the bilinear form on $U$ is positive 
definite and the assertion follows from the argument in \cite[p.~112]{Hu92}. 

{\bf Case 2:} $a^2 > 4$. In this case the bilinear form on $U$ is indefinite 
and the assertion follows immediately from Lemma~\ref{lem:a.1}.  

{\bf Case 3:} $a^2 = 4$. Here we have to distinguish two cases. 
If $\check \alpha_s$ and $\check \alpha_t$ 
are linearly dependent, then $C_S = U = \R \check \alpha_s$  
follows from $\alpha_t(\check \alpha_s) = -a < 0$ and 
$\alpha_t(\check \alpha_t) = 2 > 0$. In this case the assertion follows
from the invariance of $U$ under $\cW$. 
If $\check \alpha_s$ and $\check \alpha_t$ are linearly independent, 
then the assertion follows from Lemma~\ref{lem:a.1}. 
\end{prf}

The following proposition generalizes the corresponding 
well known result for the geometric representation 
of a Coxeter group to arbitrary linear Coxeter systems 
(cf.\ \cite[Thm.~5.4]{Hu92}). 

\begin{prop} \mlabel{prop:cox} 
Let $(V, (\alpha)_{s \in S}, (\check \alpha_s)_{s \in S})$ be a 
reflection data of finite type satisfying (C1/2). Assume that the subgroup 
$\cW \subeq \GL(V)$ generated by the reflections 
$(r_s)_{s \in S}$ is a Coxeter group. 
Let $\ell \: \cW \to \N_0$ denote its length function and 
$g \in \cW$, $s \in S$. 
\begin{description}
 \item[\rm(i)] If $\ell(gr_s) > \ell(g)$, then $g\check \alpha_s \in C_S$.
 \item[\rm(ii)] If $\ell(gr_s) < \ell(g)$, then $g\check \alpha_s \in -C_S$.
\end{description}
\end{prop}

\begin{prf} First we note that (ii) is a consequence of (i), applied to the 
element $g' = g r_s$ and using 
$g\check \alpha_s = -g'\check \alpha_s$. 

We prove (i) by induction on the length of $g$. 
The case $\ell(g) = 0$, i.e., $g = \1$ is trivial. If 
$\ell(g) > 0$, then there exists a $t \in S$ with 
$\ell(gr_t) = \ell(g)-1$. Then $t \not= s$ follows from 
$\ell(gr_s) > \ell(g)$. 
Let $\tilde S := \{ s,t\}$ and consider the subgroup 
$\tilde \cW := \la r_s, r_t \ra \subeq \cW$ with the length 
function $\tilde\ell$. Let 
\[ A := \{ f \in \cW \: f \in g\tilde \cW, \ell(g) 
= \ell(f) + \tilde\ell(f^{-1}g)  \}.\] 
Obviously $g \in A$, so that $A$ is not empty. Pick 
$f \in A$ with minimal length $\ell(f)$ and put 
$f' := f^{-1}g \in \tilde\cW$. Then $g = ff'$ with 
$\ell(g) = \ell(f) + \tilde\ell(f')$. 

We also note that $g r_t \in A$ follows from 
$\ell(g) = \ell(gr_t) + 1= \ell(gr_t) + \tilde\ell(r_t)$. 
Hence the choice of $f$ implies that 
$\ell(f) \leq \ell(gr_t) = \ell(g) -1$. We now want to 
apply the induction hypothesis to the pair 
$(f, r_s)$. To this end, we have to compare the length 
or $f$ and $fr_s$. 
If $\ell(fr_s) < \ell(f)$, then $\ell(fr_s) = \ell(f) - 1$ and we have 
\begin{align*}
 \ell(g) 
&\leq \ell(fr_s) + \ell(r_s f^{-1}g) 
\leq \ell(fr_s) + \tilde\ell(r_s f^{-1} g) \\
&\leq \ell(f) - 1+ \tilde\ell(f^{-1}g) +1 
= \ell(f) +  \tilde\ell(f^{-1}g) = \ell(g).
\end{align*}
We conclude that $\ell(g) = \ell(fr_s) + \tilde\ell(r_s f^{-1}g)$, 
so that $fr_s \in A$, contradicting the minimality of $\ell(f)$. 
This implies that $\ell(fr_s) > \ell(f)$ 
because $\ell(fr_s) \not= \ell(f)$ (\cite[p.~108]{Hu92}), 
so that the induction 
hypothesis leads to $f\check \alpha_s \in C_S$. 
By the same argument we obtain $\ell(fr_t) > \ell(f)$ and therefore 
$f\check \alpha_t \in C_S$. 

From 
\[ \ell(f) + \tilde\ell(f') = 
\ell(g) < \ell(gr_s) = \ell(ff'r_s)  
\leq \ell(f) + \ell(f'r_s)
\leq \ell(f) + \tilde\ell(f'r_s) \] 
we further derive that $\tilde\ell(f') \leq  \tilde\ell(f'r_s)$. 
Now Lemma~\ref{lem:x} implies that 
\[ f'\check \alpha_s \in C_{\tilde S} 
= \cone \{ \check \alpha_s, \check \alpha_t \}.\] 
and thus 
$g\check \alpha_s = ff'\check \alpha_s \in fC_{\tilde S} 
\subeq C_S .$
\end{prf}

\begin{thm} \mlabel{thm:d.7} 
Let $(V, (\alpha_s)_{s \in S}, 
(\check \alpha_s)_{s \in S})$ be a  linear Coxeter system 
and $\ell \: \cW \to \N_0$ be the length function with respect to the 
generating set $\{ r_s \: s \in S \}$. Then the following 
assertions hold for $g \in \cW$:  
\begin{description}
\item[\rm(i)] If $\ell(gr_s) > \ell(g)$, then $g\check \alpha_s \in C_S$ 
and $g\alpha_s \in \check C_S$. 
\item[\rm(ii)] If $\ell(gr_s) < \ell(g)$, then $g\check \alpha_s \in -C_S$ 
and $g\alpha_s \in -\check C_S$. 
\end{description}
\end{thm}

\begin{prf} In view of Proposition~\ref{prop:d.3}, the linear Coxeter system 
 $(V, (\alpha_s)_{s \in S}, (\check \alpha_s)_{s \in S})$ and its dual 
 $(V^*, (\check\alpha_s)_{s \in S}, (\alpha_s)_{s \in S})$ 
satisfy the assumptions of Proposition~\ref{prop:cox} because
the fact that $(\cW, (r_s)_{s \in S})$ is a Coxeter system implies that 
the adjoints $r_s^* \in \GL(V^*)$ define a Coxeter system in 
the subgroup $\cW \cong \la r_s^* \: s \in S \ra \subeq \GL(V^*)$. 
This implies the assertion. 
\end{prf}

\begin{rem} \mlabel{rem:pos-neg} 
For $s \in S$ and $g \in \cW$ the condition 
$g\alpha_s \in \check C_S = K^\star$ is equivalent to the linear functional 
$g\alpha_s$ taking non-negative values on $K$. 
Therefore Theorem~\ref{thm:d.7} implies in particular 
that each linear functional 
$\alpha \in \cW \{ \alpha_s \: s \in S \}$ either is 
positive or negative on~$K^0$. 
\end{rem}

\begin{prop} \mlabel{prop:d.8} 
If $(V, (\alpha_s)_{s \in S}, 
(\check \alpha_s)_{s \in S})$ is a  linear Coxeter system, $v \in K$ and 
\[ I := \{ s \in S \: \alpha_s(v) = 0\},\] 
 then the subgroup $\cW_I \subeq \cW$ generated by 
the reflections $\{r_s \: s \in I\}$ coincides with 
the stabilizer 
$\cW_v  = \{ w \in \cW \: wv = v\}.$
\end{prop}

\begin{prf} Clearly $\cW_I \subeq \cW_v$ because the generators of 
$\cW_I$  fix $v$. If $\cW_I$ is strictly smaller than $\cW_v$, 
there exists an element $g \in \cW_v \setminus \cW_I$ of minimal positive 
length. Then $\ell(g r_s) > \ell(g)$ for every $s \in I$ 
(recall that $\ell(gr_s) \not= \ell(g)$ by \cite[p.~108]{Hu92})  
implies that $g \alpha_s \in \check C_S$ (Theorem~\ref{thm:d.7}).
If $s \not\in I$, then $\alpha_s(v) > 0$ implies that 
$(g\alpha_s)(v) = \alpha_s(g^{-1}v) = \alpha_s(v) > 0$. Therefore 
$g\alpha_s$ takes positive values on $K^0$ which entails 
$g\alpha_s \in \check C_S$ (Remark~\ref{rem:pos-neg}). 
We thus arrive at $g\check C_S \subeq \check C_S$.
As the element $g^{-1} \in \cW_v \setminus \cW_I$ has the same length, 
we also obtain $g^{-1}\check C_S \subeq \check C_S$, and thus 
$g\check C_S = \check C_S$. 
Now $K = (\check C_S)^\star$ leads to $gK = K$, 
and by (LCS3) further to $g = \1$, contradicting 
$g \not\in \cW_I$. 
\end{prf}

\section{A convexity theorem for linear Coxeter systems} 
\mlabel{sec:2}

Before we come to our main theorem in this section, we have to define 
the roots and coroots of a linear Coxeter system. This is crucial to 
obtain a formulation of the theorem which does not depend on the generating 
system. This will be essential for the infinite-dimensional generalization 
where roots and coroots still make sense but $\cW$ need not be a Coxeter group.

\subsection{Roots and coroots}

\begin{defn} \mlabel{def:a.14} Let 
$(V, (\alpha_s)_{s \in S}, (\check \alpha_s)_{s \in S})$
be a linear Coxeter system. 

(a) We define the set of {\it roots} by 
\[ \Delta := \cW \{ \alpha_s \: s \in S \} \subeq V^* \quad \mbox{ and put  } \quad 
\Delta^\pm := \Delta \cap \pm K^\star 
= \Delta \cap \pm\check C_S. \] 
Roots in $\Delta^+$ are called {\it positive} and 
roots in $\Delta^-$ are called {\it negative}. 
Remark~\ref{rem:pos-neg} shows that 
\[ \Delta = \Delta^+ \dot\cup \Delta^-.\] 
We likewise define the corresponding sets of {\it coroots} 
\[ \check\Delta := \cW \{ \check \alpha_s \: s \in S\} 
\subeq V.\] 

(b) The subset 
\[ \cR := \{ w r_s w^{-1} \: w \in W, s \in S \}\subeq \cW \] 
is called the {\it set of reflections} in $\cW$. If $r = w r_s w^{-1}$ is a reflection, 
$\alpha := w \alpha_s$ and 
$\check\alpha := w \check\alpha_s$, then 
\begin{equation}
  \label{eq:genref}
r(v) = v - \alpha(v) \check \alpha \quad \mbox{ for } \quad v \in V
\end{equation}
and $\Fix(r) = \ker \alpha = w\ker \alpha_s$. 
\end{defn} 

\begin{rem}
We claim that a reflection $r \in \cW$ 
is uniquely determined by its hyperplane of fixed points. 
So let $r \in \cR$. Then the hyperplane $\Fix(r)$ intersects the interior $T^0$ of 
the Tits cone $T = \cW K$ (cf.\ Remark~\ref{rem:a.2}(a)). 
Since each $\cW$-orbit in $T^0$ meets $K$ exactly 
once (\cite[Thm.~2]{Vin71}), there exists a unique 
face $F \subeq K$ of codimension one 
and some $w \in \cW$ with 
$\Fix(w^{-1} r w) \cap K = F$. Then there exists a uniquely 
determined $s \in S$ with $F = \ker \alpha_s \cap K$ 
(Remark~\ref{rem:a.2}). 
Now $w^{-1}rw$ and $r_s$ are two reflections in the same hyperplane 
$\ker \alpha_s$, hence fixing $F$ pointwise. Therefore $w^{-1}rw(K^0) \cap r_s(K^0) \not=\eset$ 
leads to $r_s = w^{-1} r w$, so that $r = wr_s w^{-1}$. 

Next we note that, if $\alpha = w \alpha_s = w' \alpha_t$ for some 
$w,w' \in \cW$ and $s,t \in S$, then 
$w r_s w^{-1}$ and $w' r_t (w')^{-1}$ are both reflections with the 
same sets of fixed points, so that the preceding argument implies
that they are equal: $w r_s w^{-1} = w' r_t (w')^{-1}$. 
This in turn shows that 
$w\alpha_s \otimes w \check \alpha_s 
= w'\alpha_t \otimes w' \check \alpha_t 
= w\alpha_s \otimes w' \check \alpha_t$, 
and hence that $w\check \alpha_s = w' \check \alpha_t$. 
\end{rem}

\begin{defn} In view of the preceding remark, 
we can associate to each 
root $\alpha \in \Delta$ a well-defined coroot 
$\check \alpha \in \check \Delta$ such that the map 
$\Gamma \:  \Delta \to \check\Delta, \alpha \mapsto 
\check \alpha$ 
is $\cW$-equivariant and the reflections in $\cW$ have the form 
\eqref{eq:genref}. 
\end{defn}

\begin{rem}
Since the roots $\alpha_s$, $s \in S$, are 
positive by definition, Theorem~\ref{thm:d.7} implies 
that 
\[ \cone(\Delta^+) = \cone \{ \alpha_s \: s \in S \} = \check C_S,\] 
and hence that 
\[ K = (\check C_S)^\star 
= \{ v \in V \: (\forall \alpha \in \Delta^+)\, \alpha(v) \geq 0\}.\] 
\end{rem}

\subsection{Convex hulls of orbits in the Tits cone} 

The following theorem strengthens the corresponding assertion 
for elements $v \in K^0$ in \cite[p.~20]{HoG99} 
substantially because it provides also sharp information 
if the stabilizer $\cW_v$ is non-trivial. 

\begin{thm} \mlabel{thm:cox-cont} 
Let $(V, (\alpha_s)_{s \in S}, 
(\check \alpha_s)_{s \in S})$ be a  linear Coxeter system 
and $T = \cW K \subeq V$ be its Tits cone. 
For  $v \in T$ we have 
\[ \cW v \subeq v - C_v, \quad \mbox{ where } \quad 
C_v := \cone\{ \check \alpha \in\check\Delta \: \alpha(v) > 0\} .\] 
\end{thm}

\begin{prf} As $T = \cW K$ and $C_{gv} = g C_v$, 
we may w.l.o.g.\ assume that $v \in K$. 
We put \break $I := \{ s \in S \: \alpha_s(v) = 0\}$ and recall from 
Proposition~\ref{prop:d.8} that 
the corresponding parabolic subgroup $\cW_I \subeq \cW$ coincides with 
the stabilizer $\cW_v$ of~$v$. 
Let 
\[ \cW^I := \{ g \in \cW \: (\forall s \in I)\, 
\ell(gr_s) >\ell(g)\},\]  
so that $\cW = \cW^I \cW_I$ 
by \cite[p.~123]{Hu92}. 

We now show $gv - v \in - C_v$ by induction 
on the length $\ell(g)$ of $g$. 
The assertion is trivial for $g = \1$, i.e., $\ell(g) = 0$. 
Suppose that $\ell(g) > 0$. Then 
$g^{-1} = h^{-1} g_I$ with $h^{-1} \in \cW^I$ and $g_I \in \cW_I$ 
satisfying 
\[ \ell(g) = \ell(g^{-1}) = \ell(h^{-1}) + \ell(g_I) 
= \ell(h) + \ell(g_I) \] 
(\cite[p.~123]{Hu92}). 
If $g_I \not=\1$, then $\ell(h) < \ell(g)$ and 
our induction hypothesis leads to 
\[ gv - v  
= g_I^{-1} hv - v
= g_I^{-1} (hv - v) 
\in - g_I^{-1} C_v = - C_v \] 
because the stabilizer $\cW_I = \cW_v$ of $v$ 
preserves the cone $C_v$. 
We may therefore assume that 
$g_I = \1$, i.e., $g^{-1} \in \cW^I$. 
By Theorem~\ref{thm:d.7}, $g^{-1}\check C_I \subeq \check C_S$. 
In particular, $g^{-1}$ maps the set 
\[ \Delta_I^+ 
:= \Delta \cap \cone(\{ \alpha_s \: s \in I \}) 
= \{ \alpha \in \Delta^+ \: \alpha(v) = 0 \} \] 
into $\Delta^+$. 

Pick $s \in S$ with $\ell(g r_s) < \ell(g)$. 
If $\alpha_s(v) = 0$, then $gv -v = g r_s v - v \in - C_v$ 
by the induction hypothesis. We may therefore assume that 
$\alpha_s(v) > 0$.  
Then our induction hypothesis implies 
\begin{align*}
gv - v 
&= (gr_s)r_s v - v 
= (gr_s)(v - \alpha_s(v) \check \alpha_s) - v \\
&= (gr_s)v - v - \alpha_s(v) (gr_s)\check \alpha_s 
\in - C_v + \alpha_s(v) g\check \alpha_s. 
\end{align*}
In view of $\alpha_s(v) > 0$, it remains to see that 
$g\check \alpha_s = (g\alpha_s)\,\check{} \in -C_v$ 
(cf.\ Definition~\ref{def:a.14}(c)), so that it suffices 
to verify that $(g\alpha_s)(v) < 0$. 
As $g(-\alpha_s) \in \Delta^+ \subeq \check C_S$ by 
Theorem~\ref{thm:d.7}, we have $(g\alpha_s)(v) \leq 0$. 
If $(g\alpha_s)(v) = 0$, then 
\[ -\alpha_s \in 
g^{-1}\{ \beta \in \Delta^+ \: \beta(v) =0\} 
= g^{-1}\Delta^+_I \subeq \Delta^+\] 
by construction of $g$; but this leads to the contradiction 
$-\alpha_s \in \Delta^+$. This proves that 
$(g\alpha_s)(v) < 0$, which completes the proof.
\end{prf}

\begin{cor} \mlabel{cor:1.16} For $v \in T$, the following assertions hold: 
  \begin{description}
  \item[\rm(i)] $\cone(\cW v - v) = - C_v.$ 
  \item[\rm(ii)]  $v$ is an extreme point of $\co(v)$ if and only if 
the cone $C_v$ is pointed. 
\item[\rm(iii)] For $\lambda \in V^*$, 
\[ \lambda(v) = \min\lambda(\cW v) 
\quad \Leftrightarrow \quad \lambda \in -C_v^\star.\] 
  \end{description}
\end{cor}

\begin{prf} (i) For any $v \in V$ with 
$\alpha(v) > 0$ the relation 
$r_\alpha(v) = v - \alpha(v) \check \alpha$ implies that 
\[ - \check \alpha \in \R_+ (\cW v - v),\] 
so that Theorem~\ref{thm:cox-cont} implies for $v \in T$ that 
$\cone(\cW v - v) = - C_v.$ 

(ii) and (iii) follow immediately from (i). 
\end{prf}

\begin{cor} \mlabel{cor:1.17} The following conditions are equivalent 
  \begin{description}
    \item[\rm(i)] The cone $C_S$ is pointed, i.e., 
$(V^*, (\check\alpha_s)_{s \in S}, (\alpha_s)_{s \in S})$ 
is a linear Coxeter system. 
    \item[\rm(ii)] There exists a $v \in K^0$ which is an extreme point 
of $\co(v)$. 
    \item[\rm(iii)] Each  $v \in K^0$ is an extreme point of $\co(v)$. 
  \end{description}
\end{cor}

\begin{prf} This is immediate from Corollary~\ref{cor:1.16}  and 
the fact that $C_S = C_v$ for $v \in K^0$. 
\end{prf}

\begin{prob} From Theorem~\ref{thm:cox-cont} 
we obtain for $v \in T$ the relation 
\[ \co(v) \subeq \bigcap_{w \in \cW} w(v - C_v).\] 
When do we have equality? 
\end{prob}

\begin{ex} We take a closer look at linear Coxeter systems with a 
$2$-element set $S = \{s,t\}$ and $\dim V = 2$. 

First we show that $\alpha_s$ and $\alpha_t$ are linearly independent. 
If this is not the case, then $\alpha_s(\check \alpha_t) \leq 0$ implies that 
$\alpha_s = \lambda \alpha_t$ for some $\lambda < 0$, but this leads to the 
contradiction  $K = \eset$. Therefore $\alpha_s$ and $\alpha_t$ 
are linearly independent. 

(a) Suppose that 
$\alpha_s(\check \alpha_t) = \alpha_t(\check \alpha_s) = -2$. 
Then $\alpha_s$ and $\alpha_t$ vanish on $\check\alpha_s + \check \alpha_t$,  
and since $V^*$ is spanned by $\alpha_s$ and $\alpha_t$, it follows 
that $\check\alpha_t = - \check \alpha_s$. In particular, the cone 
$C_S = \R \check \alpha_s$ is not pointed. 

This implies that the action of $\cW$ on $V$ leaves all 
affine subspaces of the form $v + C_S$ invariant. 
If $\alpha_s^*, \alpha_t^* \in V$ is the dual basis of 
$\alpha_s, \alpha_t$, then 
\[ K = \R_+ \alpha_s^* + \R_+ \alpha_t^* \quad \mbox{ and } \quad 
\check \alpha_s = - \check \alpha_t = 2(\alpha_s^* - \alpha_t^*).\] 
The linear map $r_t r_s$ fixes the line $C_S$ pointwise and induces on 
$V/C_S$ the identity, hence is unipotent. Moreover, 
\[ r_t r_s(\alpha_s^*) 
= r_t(\alpha_s^* - \check \alpha_s)
= r_t(-\alpha_s^* + 2 \alpha_t^*) 
= -\alpha_s^* + 2 \alpha_t^* - 2 \check \alpha_t
= -5\alpha_s^* + 6 \alpha_t^* 
= \alpha_s^* - 3 \check \alpha_s, \] 
so that the convexity of the Tits cone implies that 
\[ T^0 = C_S + \R_+^\times \alpha_s^* \] 
is an open half plane and $T = T^0 \cup \{0\}$. We further have 
\[ \co(v) = v + C_S \quad \mbox{ for } \quad v \in T^0\] 
and 
\[ \check \Delta = \{ \check \alpha_s, \check \alpha_t\}  
= \{ \pm \check \alpha_s\}, \] 
whereas $\Delta = \cW \{ \alpha_s, \alpha_t\}$ is infinite. 
More precisely, we have 
\[ r_t r_s \alpha_s 
= r_t(-\alpha_s) 
= - \alpha_s + \alpha_s(\check \alpha_t) \alpha_t 
= - \alpha_s - 2 \alpha_t 
= \alpha_s - 2(\alpha_s + \alpha_t), \] 
and since $\cW = \{ r_s (r_t r_s)^n, (r_t r_s)^n \: n \in \Z \}$, 
it follows that 
\[ \Delta = \{ \pm \alpha_s, \pm \alpha_t \} 
+ 2 \Z (\alpha_s + \alpha_t).\] 
The two coroots $\check \alpha_s$ and $\check \alpha_t$ lie in the  
boundary of the Tits cone and 
$\cW \check \alpha_s = \{ \pm \check \alpha_s\}.$ 
Moreover, 
\[ \{ \alpha \in \Delta \: \alpha(\check \alpha_s) > 0 \} 
= \{ \alpha_s, - \alpha_t \} + 2\Z(\alpha_s + \alpha_n), \] 
so that 
\[ C_{\check \alpha_s} 
= \cone\{ \check \alpha \: \alpha(\check \alpha_s) > 0\} 
= \R_+ \check \alpha_s \] 
has the property that 
\[ \cW \check \alpha_s \subeq \check \alpha_s - C_{\check \alpha_s}.\] 
That this remains true on the boundary $\partial T$ is due to the fact that 
the $\cW$-action on this line can also be obtained from the 
one-dimensional reflection  datum $(\partial T, \alpha_s, \check \alpha_s)$. 

(b) Suppose that 
$\alpha_s(\check \alpha_t) = \alpha_t(\check \alpha_s) = -3$. 
Then 
\[ \det\pmat{ 
\alpha_s(\check \alpha_s) & \alpha_s(\check \alpha_t) \\
\alpha_t(\check \alpha_s) & \alpha_t(\check \alpha_t)} 
= 4 - 9 < 0 \] 
implies that $\check \alpha_s$ and $\check \alpha_t$ are linearly 
independent. 
In this case the symmetric bilinear form $(\cdot,\cdot)$ 
on $V$ represented by the matrix 
\[ A = \pmat{ 1 & -\frac{3}{2} \\ -\frac{3}{2} & 1} \] 
with respect to the basis $\check \alpha_s, \check \alpha_t$ 
is $\cW$-invariant. 
If $\alpha_s^*, \alpha_t^* \in V$ is the dual basis of 
$\alpha_s, \alpha_t$, then 
\[ K = \R_+ \alpha_s^* + \R_+ \alpha_t^* \quad \mbox{ and } \quad 
\alpha_s^* = -\frac{2}{5} \check \alpha_s - \frac{3}{5}\check \alpha_t, \quad 
\alpha_t^* = -\frac{3}{5} \check \alpha_s - \frac{2}{5}\check \alpha_t.\] 
Now $(\alpha_s^*, \alpha_s^*) =(\alpha_t^*, \alpha_t^*) < 0$ implies that 
\[ T = \cW K \subeq 
\{ v \in V \: (v,v) \leq 0, (v, \alpha_s^*) \leq 0 \}. \] 
In particular, $\check \alpha_s$ is not contained in 
$\pm T$. For this element we have 
\[ r_s r_t(\check\alpha_s) 
= r_s(\check \alpha_s + 3 \check \alpha_t) 
= - \check \alpha_s + 3 (\check \alpha_t + 3 \check \alpha_s) 
=  8\check \alpha_s + 3 \check \alpha_t \] 
and 
\[ r_t r_s(\check\alpha_s) 
= r_t(-\check \alpha_s) 
= - \check \alpha_s - 3 \check \alpha_t, \] 
so that 
\[ \frac{7}{2}\check \alpha_s  
= \frac{1}{2}(r_t r_s(\check\alpha_s)  + r_s r_t(\check\alpha_s)) 
\in \co(\check \alpha_s) \] 
shows that $\check \alpha_s$ is an interior point of 
$\co(\check \alpha_s)$. 

On the other hand, the cone 
\[ C_{\check \alpha_s} 
= \cone\{ \check \alpha \: \alpha(\check \alpha_s) > 0 \}
= \cone\{ \check \alpha \: (\check \alpha, \check \alpha_s) > 0 \} \] 
is proper, so that 
$\cW\check \alpha_s \not\subeq \check \alpha_s - C_{\check \alpha_s}.$ 
\end{ex}

In general the dual of a linear Coxeter system is not a linear Coxeter system. 
However, we have seen in Remark~\ref{rem:1.7}(b) that we always obtain a linear Coxeter system 
on the subspace $U = \Spann(C_S^\star) \subeq V^*$. For orbits in the corresponding Tits cone, 
we have the following variant of Theorem~\ref{thm:cox-cont}. 

\begin{thm} \mlabel{thm:cox-cont2} 
{\rm(Convexity Theorem for $V^*$)} 
Let $(V, (\alpha_s)_{s \in S}, 
(\check \alpha_s)_{s \in S})$ be a  linear Coxeter system  
and $\check T = \cW C_S^\star \subeq V^*$. 
For any $\lambda \in \check T$ we then have 
\[ \cW \lambda \subeq \lambda  - C_\lambda \quad \mbox{ for } \quad 
C_\lambda := \cone\{ \alpha \in\Delta \: \lambda(\check \alpha) > 0\} .\] 
\end{thm}

\begin{prf} We may w.l.o.g.\ assume that 
$\lambda \in C_S^\star$. Consider the subspace 
\[ U := \Spann(C_S^\star) = H(C_S)^\bot \subeq V^*\] 
and recall from Remark~\ref{rem:1.7}(b) that 
$(U, (q(\check \alpha_s))_{s \in \tilde S}, (\alpha_s)_{s \in \tilde S})$ 
is a linear Coxeter system with fundamental chamber 
$C_S^\star$ for which we have the surjective restriction map 
\[ R \: \cW \to \cW_U, \quad w\mapsto w\res_U. \] 
It follows in particular that, 
for $\lambda \in C_S^\star \subeq U$, we have 
$\cW \lambda = \cW_U \lambda.$ 
Applying Theorem~\ref{thm:cox-cont}, we obtain 
\begin{equation}\label{eq:wrel}
\cW\lambda = \cW_U \lambda 
\subeq \lambda - C_\lambda^U,\mbox{ where } \quad 
C_\lambda^U = \cone\{ \alpha \in \check\Delta_U \: \lambda(\check\alpha)>0\}.
\end{equation}
On the other hand $\cW\res_U =\cW_U$ implies 
\[ \check\Delta_U 
= \cW_U \{ \alpha_s \: s \in \tilde S\} 
= \cW \{ \alpha_s \: s \in \tilde S\} \subeq U.\] 
For $s\in S \setminus \tilde S$ we have 
$\check \alpha_s \in U^\bot$, hence also 
$\cW \check \alpha_s \in U^\bot$ because $U$ is $\cW$-invariant 
and thus $\lambda(\cW\check \alpha_s) =\{0\}$. 
This shows that 
\[ C_\lambda = \cone\{ \alpha \in \check\Delta \: \lambda(\check\alpha)>0\} 
= C_\lambda^U,\] 
and by \eqref{eq:wrel}, the proof is complete. 
\end{prf}

The following proposition extends Proposition~\ref{prop:d.8} 
to stabilizers of elements in $C_S^\star$. 

\begin{prop} \mlabel{prop:1.23} If $\lambda \in C_S^\star$, then 
$\cW_\lambda = \la r_s \: \lambda(\check \alpha_s) = 0 \ra.$
\end{prop}

\begin{prf} We recall the subspace 
$U := (\check C_S)^\star - (\check C_S)^\star \subeq V^*$ 
and the related objects also used in the preceding proof. 
We then have a surjective homomorphism
$R \: \cW \to \cW_U, w \mapsto w\res_U,$ 
and $\lambda \in U$ implies that $\ker R \subeq \cW_\lambda$. 

For $s \in S \setminus \tilde S$ we have 
$\check \alpha_s \in H(C_S)\subeq \ker \lambda$, 
which leads to $r_s \in \ker R \subeq \cW_\lambda$. Therefore 
\[ S_\lambda := \{ s \in S \: \lambda(\check \alpha_s) = 0\} 
\supeq S \setminus \tilde S.\] 
Let $\tilde S_\lambda :=\tilde S \cap S_\lambda$. 
Since $C_S^\star$ is the fundamental chamber of the 
linear Coxeter system in $U$, Proposition~\ref{prop:d.8} yields 
\[ \cW_{U,\lambda} 
= \la r_s \: s \in \tilde S, \lambda(\check \alpha_s) =0\ra.\] 
This implies that 
$\cW_\lambda \subeq \la r_s \: s \in \tilde S_\lambda\ra \cdot \ker R.$ 

Next we observe that, for $s \in \tilde S$ and 
$t \in S \setminus \tilde S$ we have 
$\alpha_s \in U$ and $\check \alpha_t \in H(C_S) = U^\bot$, so that 
$\alpha_s(\check \alpha_t) =0$. From (C1) we now also obtain 
$\alpha_t(\check \alpha_s) =0$ and this implies that 
$r_s r_t = r_t r_s$: 
\begin{align*}
r_s r_t(v) 
&= r_s(v - \alpha_t(v)\check \alpha_t)
= v - \alpha_t(v)\check \alpha_t 
- \alpha_s(v) \check \alpha_s 
+ \alpha_t(v) \alpha_s(\check \alpha_t) \check \alpha_s\\
&= v - \alpha_t(v)\check \alpha_t 
- \alpha_s(v) \check \alpha_s 
+ \alpha_s(v) \alpha_t(\check \alpha_s) \check \alpha_t= r_t r_s(v).
\end{align*}
Therefore 
\[ \cW 
= \la r_s \: s \in S \ra 
= \la r_s \: s \in \tilde S \ra \la r_t \: t \in S \setminus \tilde S \ra
= \cW_{\tilde S} \cW_{S \setminus \tilde S}\]  
is a product of two commuting subgroups. 
Since the subgroup $\cW_{\tilde S}$ of $\cW$ is a Coxeter group 
with Coxeter system $\{ r_s \: s \in \tilde S\}$, 
the restriction homomorphism $R \: \cW \to \cW_U$ maps 
$\cW_{\tilde S}$ bijectively onto $\cW_U$. 
On the other hand, $\cW_{S\setminus \tilde S} \subeq \ker R$, so that 
$\ker R \cap \cW_{\tilde S} = \{\1\}$ leads to 
$\ker R = \cW_{S\setminus \tilde S}.$ 
We finally arrive at 
\[ \cW_\lambda 
= \cW_{\tilde S_\lambda} \ker R 
= \cW_{\tilde S_\lambda} \cW_{S \setminus \tilde S} 
= \cW_{S_\lambda}.\qedhere\] 
\end{prf}

The following proposition describes the subset of 
a $\cW$-orbit in $T$ on which a linear functional 
$\lambda \in C_S^\star$ takes its maximal values as the orbit 
of the stabilizer of $\cW_\lambda$ and we also provide a dual version. 

\begin{prop}\mlabel{prop:1.24}
Let $\lambda \in C_S^\star$ and $v \in K = (\check C_S)^\star$, 
so that 
\[ \lambda(v) = \max \lambda(\cW v).\] 
If $g \in \cW$ satisfies $\lambda(gv) = \max \lambda(\cW v)$, 
then $gv \in \cW_\lambda v$ and $g^{-1}\lambda \in \cW_v \lambda$. 
\end{prop}

\begin{prf} (a) First we show that $gv \in \cW_\lambda v$. 
In Proposition~\ref{prop:1.23} we have seen that 
$\cW_\lambda = \cW_{S_\lambda}$ is a parabolic subgroup of $\cW$. 
Let 
\[ \cW^\lambda := \{ w \in \cW \: (\forall s \in S_\lambda)\, 
\ell(r_sw) > \ell(w)\},\]  
so that $\cW = \cW_\lambda \cW^\lambda$ by \cite[p.~123]{Hu92}. 

If there exists a $g \in \cW$ with $gv \not\in \cW_\lambda v$ and 
$\lambda(gv) = \max \lambda(\cW v)$, we choose  such an element 
$g$ with minimal length. Then 
$g \in \cW^\lambda$ with $\ell(g) > 0$.
We pick an $s \in S$ with 
$\ell(g^{-1}r_s) = \ell(r_sg) < \ell(g)$ and observe that 
this implies that $s \not\in S_\lambda$, i.e., 
$\lambda(\check \alpha_s) \not=0$ and therefore 
$\lambda(\check \alpha_s) > 0$ because $\lambda \in C_S^\star$. 
Then $g^{-1}\alpha_s \in -\check C_S$ by 
Theorem~\ref{thm:d.7}, which leads to 
$0 \geq (g^{-1}\alpha_s)(v) = \alpha_s(gv)$. We 
thus arrive at 
\[ \lambda(r_s gv) 
= \lambda(gv) - \underbrace{\alpha_s(gv)}_{\leq 0} 
\underbrace{\lambda(\check \alpha_s)}_{> 0} 
\geq \lambda(gv) \geq \lambda(r_sgv),\] 
where the last inequality follows from the maximality of 
$\lambda(gv)$. We conclude that 
$\alpha_s(gv) =0$, so that $r_sgv = gv \not\in \cW_\lambda v$. 
This contradicts the minimality of the length of $g$. 

(b) Now we show that $g^{-1}v \in \cW_v \lambda$. 
In Proposition~\ref{prop:d.8} we have seen that 
$\cW_v$ is a parabolic subgroup of $\cW$ generated by the 
reflections $r_{\alpha_s}$ with 
\[ s\in S_v := \{ s\in S \: \alpha_s(v) = 0\}.\] 
Let 
\[ \cW^v := \{ w \in \cW \: (\forall s \in S_v)\, 
\ell(wr_s) > \ell(w)\},\]  
so that $\cW = \cW^v\cW_v$ by \cite[p.~123]{Hu92}. 

If there exists a $g \in \cW$ with $g^{-1}\lambda \not\in \cW_v\lambda$ and 
$\lambda(gv) = \max \lambda(\cW v)$, we choose  such an element 
$g\in \cW$ of minimal length. Then 
$g \in \cW^v$ with $\ell(g) > 0$.
We pick an $s \in S$ with 
$\ell(gr_s)  < \ell(g)$ and observe that 
this implies that $s \not\in S_v$, i.e., 
$\alpha_s(v) \not=0$ and therefore 
$\alpha_s(v) > 0$ because $v \in K$. 
We further obtain $g\check\alpha_s \in -C_S$ from 
Proposition~\ref{prop:cox}, which leads to 
$0 \geq \lambda(g\check\alpha_s) = (g^{-1}\lambda)(\check \alpha_s).$ We 
thus arrive at 
\[ \lambda(gr_sv ) 
= \lambda(gv) - \underbrace{\alpha_s(v)}_{> 0} 
\underbrace{\lambda(g\check \alpha_s)}_{\leq  0} 
\geq \lambda(gv) \geq \lambda(gr_sv),\] 
where the last inequality follows from the maximality of 
$\lambda(gv)$. We conclude that 
$\lambda(g\check \alpha_s) =0$, so that 
$r_s(g^{-1}\lambda) = g^{-1}\lambda \not\in \cW_v\lambda$. 
This contradicts the minimality of the length of $g$. 
\end{prf}

\begin{prob} For $v \in K$ and $\alpha \in \Delta$, the relation 
$\alpha(v) > 0$ implies $\alpha \in \Delta^+ \subeq \check C_S$, 
so that $C_v \subeq C_S$ and 
$C_S^\star \subeq C_v^\star$. Is the conclusion of 
Proposition~\ref{prop:1.24} still valid under the weaker  
assumption $\lambda \in C_v^\star$? 
\end{prob}

\section{Orbits of locally finite and locally affine Weyl groups} 
\mlabel{sec:3}

Now we turn to the applications to Weyl group orbits 
of linear functionals for locally finite and locally affine root systems 
(Theorems~\ref{thm:cox-cont} and \ref{thm:cox-cont2}). 

\subsection{Locally finite  root systems} 
\mlabel{subsec:3.1}

First we describe the irreducible locally 
finite root systems of infinite rank (cf.\ \cite[\S 8]{LN04}, 
\cite{NS01}). 
Let $J$ be a set and $V := \R^{(J)}$ denote the free vector space over 
$J$, endowed with the canonical scalar product, given by 
\[ (x,y) := \sum_{j \in J} x_j y_j.\] 
We write $(e_j)_{j \in J}$ for the canonical orthonormal basis and we 
realize the root systems in the dual space $V^* \cong \R^J$ which contains
the linearly independent system $\eps_j := e_j^*$, defined by 
$\eps_j^*(e_k) = \delta_{jk}$. On $\Spann \{ \eps_j \: j \in J\}$ we also have a positive 
definite scalar product defined by $(\eps_i, \eps_j) = \delta_{ij}$ 
for which the canonical inclusion $V \into V^*$ is isometric. 
The infinite irreducible locally finite root systems are given by 
\begin{align*}
A_{J} &:= \{ \eps_{j} - \eps_{k} : j,k \in J, j \not= k \}, \\
B_{J} &:= \{ \pm \eps_{j}, \pm \eps_{j} \pm \eps_{k} : j,k \in J, j \not= k \},\\
C_{J} 
&:= \{ \pm 2 \eps_{j}, 
\pm \eps_{j} \pm \eps_{k} : j,k \in J, j \not= k \},\\ 
D_{J} &:= \{ \pm \eps_{j} \pm \eps_{k} : j,k \in J, j \not= k \}, \\ 
BC_{J} &:= \{ \pm \eps_{j}, \pm 2\eps_{j}, 
\pm \eps_{j} \pm \eps_{k} : j,k \in J, j \not= k \}. 
\end{align*}

Let $\Delta \subeq V^* \cong \R^J$ be a locally finite root system 
of type $X_J$ with $X \in \{A,B,C,D,BC\}$. 
For $\alpha \in \Spann \Delta$, we write 
$\alpha^\sharp \in V$ for the unique element determined by 
\[ \alpha(v) = (v, \alpha^\sharp) \quad \mbox{ for } \quad v \in V.\] 
For $\alpha \in \Delta$ we define its {\it coroot} by 
\[ \check \alpha = \frac{2}{(\alpha,\alpha)} \alpha^\sharp.\] 
This leads to a reflection system $(V,\Delta, \check\Delta)$ and a corresponding 
group $\cW$, called in this context the {\it Weyl group}. 

\begin{thm} \mlabel{thm:conv-locfin} For $\lambda \in V^*$ we have 
\[ \cW \lambda \subeq \lambda  - C_\lambda \quad \mbox{ for } \quad 
C_\lambda := \cone\{ \alpha \in\Delta \: \lambda(\check \alpha) > 0\} .\] 
\end{thm}

\begin{prf} Let $w \in \cW$ and observe that $w$ is a finite product 
of reflections $r_{\alpha_1}, \ldots, r_{\alpha_n}$. We pick a finite 
dimensional subset $F \subeq J$ such that 
$\alpha_j \in \R^F$ for $j =1,\ldots, n$. Accordingly, we have an orthogonal 
direct sum 
\[ V = V_0 \oplus V_1 \quad \mbox{ with } \quad 
V_0 = \R^F \quad \mbox{ and } \quad V_1 := V_0^\bot = \R^{J \setminus F}\] 
which is invariant under $w$. 

Next we observe that $\Delta_0 := \Delta \cap \R^F$ is a finite root system 
of type $X_F$, which implies that the finite reflection 
system $(V_0, \Delta_0, \check \Delta_0)$ comes from a finite Coxeter system 
with finite Coxeter group $\cW_0$ containing $w$. 
In this case the Tits cones in $V_0$ and $V_0^*$ coincide with the whole 
space, so that 
\[ \cW_0 \lambda_0 \subeq \lambda_0 - C_{\lambda_0}^0 \] 
holds for $\lambda_0 := \lambda\res_{V_0}$ and 
\[ C_{\lambda_0}^0 
:= \cone\{ \alpha \in\Delta_0 \: \lambda_0(\check \alpha) = \lambda(\check \alpha)> 0\} 
\subeq C_\lambda.\] 
Writing $\lambda  = \lambda_0 \oplus \lambda_1$ according to the decomposition 
$V = V_0 \oplus V_1$, we now obtain 
\[ w\lambda = w\lambda_0 \oplus \lambda_1 
\in (\lambda_0- C_{\lambda_0}^0) \oplus \lambda_1
= \lambda - C_{\lambda_0}^0
\subeq \lambda - C_\lambda.\qedhere\] 
\end{prf}

\begin{cor} For $d \in V$ and $\lambda \in V^*$, the following are equivalent 
  \begin{description}
    \item[\rm(i)] $\lambda(d) = \min \la \cW\lambda, d\ra$. 
    \item[\rm(ii)] $d \in - C_\lambda^\star$. 
    \item[\rm(iii)] $(\forall \alpha \in \Delta)\ \lambda(\check \alpha) > 0 
\Rarrow \alpha(d) \leq 0$. 
  \end{description}
\end{cor}

\begin{rem} Since the canonical inclusion 
$V = \R^{(J)} \into V^*\cong \R^J$ is $\cW$-equivariant, Theorem~\ref{thm:conv-locfin}  
implies the corresponding result for $\cW$-orbits in $V$ itself: 
\[ \cW v \subeq v - C_v \quad \mbox{ for } \quad 
C_v = \cone \{ \check \alpha \: \alpha(v) > 0 \}.\]  
\end{rem}

\subsection{Locally affine root systems} 
\mlabel{subsec:3.2} 

Let $V = \R^{(J)}$ be as above and $\Delta \subeq V^* \cong \R^J$ be a 
locally finite root system of type $X_J$. 
We put 
\[ \hat V := \R \times V \times \R \quad \mbox{ and } \quad 
\Delta^{(1)} := \{0\} \times \Delta \times \Z  
\subeq \R \times V^* \times \R \cong \hat V^*. \] 
We also define a Lorentzian form on $\hat V$ by 
\[ ((z,x,t), (z,x', t')) := (x,x') - zt' - z't.\] 
Suppressing the first component, we have Yoshii's classification (\cite[Cor.~13]{YY10}): 

\begin{prop} \mlabel{prop:classlocaff} 
The irreducible reduced locally affine root systems of infinite rank 
are the following, where $J$ is an infinite set: 
$A_{J}^{(1)}, B_J^{(1)}, C_J^{(1)}, D_J^{(1)},$ 
or 
\begin{align*}
B_J^{(2)} &:= \big(B_J\times 2\Z\big) 
\cup \big(\{ \pm \eps_j \: j \in J\}\times (2\Z + 1)\big), \\
C_J^{(2)} &:= (C_J \times 2\Z) 
\cup \big(D_J  \times (2 \Z+1)\big) \\
(BC_J)^{(2)} &:= (B_J \times 2\Z) \cup \big(BC_J\times (2\Z+1)\big). 
\end{align*}
\end{prop}

Let $\hat\Delta \subeq \hat V^*$ be one of these locally affine root systems. 
We write 
\[ \Delta_n := \{ \alpha \in \Delta \: (0,\alpha,n) \in \hat\Delta \}, \] 
so that 
\[ \hat\Delta = (\{0\} \times \Delta_0 \times 2 \Z) 
\dot\cup (\{0\} \times \Delta_1 \times (2\Z + 1)).\] 
A quick inspection shows that all reflections corresponding to 
roots in $\Delta_1$ are also obtained from $\Delta_0$. Therefore we obtain an injection 
\[ \iota_\cW \: \cW \cong \la r_{(0,\alpha,0)} 
\: \alpha \in \Delta_0\ra_{\rm grp} \into \hat \cW.\] 

Since $\hat\Delta$ consists of non-isotropic vectors for the Lorentzian form, we can also 
define for $\uline\alpha = (0, \alpha,n) \in \hat\Delta$ the {\it coroot} by 
\begin{equation}
  \label{eq:coroot}
 \check{\uline\alpha} 
= \frac{2}{(\uline\alpha,\uline\alpha)} {\uline\alpha}^\sharp
= \frac{2}{(\alpha,\alpha)}(-n, \alpha^\sharp,0)
= \Big(\frac{-2n}{(\alpha,\alpha)}, \check \alpha,0\Big)
\end{equation}
and obtain a reflection system $(\hat V,\hat\Delta, 
\check{\hat\Delta})$ and a corresponding (affine) Weyl group~$\hat\cW$. 
In the following we write 
\[ c := (1,0,0) \quad \mbox{ and } \quad d := (0,0,1)\]  
for these two distinguished elements of $\hat V$. 

\begin{thm} \mlabel{thm:convlocaff} 
For $\lambda\in \hat V^*$ with $\lambda(c)\not=0$ the following assertions hold: 
\begin{description}
\item[\rm(i)] $\hat\cW \lambda \subeq \lambda  - C_\lambda$ for 
$C_\lambda := \cone\{ \alpha \in\hat\Delta \: \lambda(\check \alpha) > 0\}$. 
\item[\rm(ii)] If $\lambda(d) = \min (\hat\cW \lambda)(d)$, 
then any $\mu \in \hat\cW\lambda$ with $\mu(d) = \lambda(d)$ is contained 
in the orbit of $\cW \cong \hat\cW_d$.
\end{description}
\end{thm}

\begin{prf} (i) Since we can argue as in the locally finite case, it suffices to 
show that, if $J$ is finite and $\hat\Delta$ of type $X_J^{(t)}$ for $t = 1,2$, then 
$\lambda \in \hat V^*$ is contained in the Tits cone, so that 
Theorem~\ref{thm:cox-cont2} applies. To this end, we have to recall 
the description of affine root systems with respect to a simple 
system of roots. 

Since $\lambda(c) \not=0$, the $\hat\cW$-orbit contains an element which 
is either dominant or antidominant 
(\cite[Prop.~4.9]{Ne10}; see also \cite[Prop.~6.6]{Ka90} and 
\cite[Thm.~16]{MP95}), i.e., 
there exists a simple systems of roots $\Pi \subeq \hat\Delta$ 
such that $\lambda$ is dominant. This means that 
$\lambda$ is contained in the corresponding fundamental chamber, hence in 
particular in the Tits cone. Now the assertion follows from 
Theorem~\ref{thm:cox-cont}. 

(ii) With the same argument as before, it suffices it suffices to 
prove the assertion for the case where $J$ is finite and 
$\lambda$ is antidominant with respect to a given simple system $\Pi$ of roots. 
So assume that $\hat\Delta$ is of type $X_J^{(t)}$ 
for some finite set $J$. 
We write 
$\Pi = \{ \alpha_0, \ldots, \alpha_r \}$ 
for a set of simple roots, where 
$\alpha_1, \ldots, \alpha_r$ are simple roots of the corresponding 
finite root system $\Delta$ and 
$\alpha_0 = (-\theta, 1)$, where 
$\theta$ is the ``highest weight'' in $\Delta_1$ 
with respect to the positive system defined by $\{ \alpha_1, \ldots, \alpha_r\}$ 
(cf.\ \cite{Ka90}). 

Now 
$\alpha_j(d) = 0$ for $j = 1,\ldots, r$ and 
$\alpha_0(d) = 1$ imply that $d \in K = (\check C_S)^\star$. 
Further, the antidominance of $\lambda$ with respect to $\Pi$ means that 
$\lambda \in -C_S^\star$. 
If  $\mu = w^{-1}\lambda$ satisfies 
$\mu(d) = \lambda(wd) = \lambda(d) = \min (\hat\cW\lambda)(d)$, 
we obtain $\mu  \in \hat\cW_d \lambda$ from Proposition~\ref{prop:1.24}. 

Finally, we note that the stabilizer group $\hat\cW_d$ 
is a parabolic subgroup of $\hat\cW$ generated by 
the fundamental reflections $r_{\alpha_j}$ fixing $d$, which is the case 
for $j > 0$. Therefore $\hat\cW_d \cong  \cW$ is the Weyl group of the 
corresponding finite root system $X_J$. 
\end{prf}

\begin{cor} \mlabel{cor:2} 
For $d = (0,0,1)$ and $\lambda \in \hat V^*$, the following are equivalent 
  \begin{description}
    \item[\rm(i)] $\lambda(d) = \min \la \hat\cW\lambda, d\ra$. 
    \item[\rm(ii)] $d \in - C_\lambda^\star$. 
    \item[\rm(iii)] $(\forall \uline\alpha \in \hat\Delta)\ \lambda(\check{\uline\alpha}) > 0 
\Rarrow \uline\alpha(d) \leq 0$. 
    \item[\rm(iv)] $(\forall \uline\alpha = (0,\alpha,n)\in \hat\Delta)\ 
n > 0 \Rarrow  \frac{(\alpha,\alpha)}{2n}\lambda(\check \alpha) \leq \lambda(c).$ 
  \end{description}
\end{cor}

\begin{prf} The equivalence of (i) and (ii) follows from 
Theorem~\ref{thm:convlocaff}(i), and (iii) is a reformulation  
of (ii). The equivalence of (iii) and (iv) follows by negating the implication, 
inserting the formula for the coroot and using $\uline\alpha(d) = n$. 
\end{prf}

\begin{defn} Linear functionals $\lambda \in \hat V^*$ satisfying the equivalent 
conditions in Corollary~\ref{cor:2} are called {\it $d$-minimal}. 
\end{defn}

\subsection{Characterization of $d$-minimal weights} 
\mlabel{subsec:3.3} 

For $\lambda \in \hat V^*$, let $\lambda_c := \lambda(c)$. 

\begin{lem} \mlabel{lem:4.2} If $(\hat\cW\lambda)(d)$ is bounded from 
below, then  $\lambda_c \geq 0$. If, in addition, 
$\lambda_c = 0$, then $\lambda$ is fixed by $\hat\cW$. 
\end{lem}

\begin{prf} If $(0,\alpha,n) \in \hat\Delta$, then also 
$(0,\alpha,n + 2k) \in \hat\Delta$ for every $k\in \N$, so that 
Corollary~\ref{cor:2}(iv) implies $\lambda_c \geq 0$.

If $\lambda_c = 0$, 
then $(0.\alpha,n) \in\hat\Delta$ for some $n \not=0$ implies 
the additional condition $\lambda(\check \alpha) =0$. 
Hence 
$\lambda_c = 0$ leads to $\lambda(\check{\uline\alpha}) = 0$ for every 
$\uline\alpha \in \hat\Delta$, so that $\lambda$ is fixed by 
$\hat\cW$.   
\end{prf}

\begin{prop} \mlabel{prop:4.2} 
Suppose that $\hat\Delta$ is one of the $7$ irreducible 
locally affine root systems. 
For $\lambda \in V^*$ with 
$\lambda_c> 0$, the following are equivalent: 
\begin{description}
\item[\rm(i)] $\lambda$ is $d$-minimal. 
\item[\rm(ii)] $(\forall \alpha \in \Delta, n =1,2)\quad 
(0,\alpha,n) \in\hat\Delta \Rarrow |\lambda(\check\alpha)| \frac{(\alpha,\alpha)}{2n}\leq \lambda_c$. 
\end{description}
\end{prop}

\begin{prf} That (ii) follows from (i) is a consequence of Corollary~\ref{cor:2}(iv) 
and the observation that $(0,\alpha,n) \in \hat\Delta$ implies 
$(0,-\alpha,n) \in \hat\Delta$. 

If, conversely, (ii) holds, then the $2$-periodic structure 
of the root system implies the condition in Corollary~\ref{cor:2}(iv). 
\end{prf}

\begin{thm} \mlabel{thm:dmin} 
For the seven irreducible locally affine 
root systems $X_J^{(r)} = \hat\Delta$ of infinite rank, 
a linear functional $\lambda = (\lambda_c, \oline \lambda, \lambda_d) 
\in \hat V^*$ with $\lambda_c > 0$ is $d$-minimal 
if and only if the following conditions are satisfied by 
the corresponding function $\oline\lambda \: J \to \R, j \mapsto \lambda_j$: 
\begin{description}
\item[\rm($A_J^{(1)}$)] 
$\max\oline\lambda - \min \oline\lambda \leq\lambda_c$. 
\item[\rm($B_J^{(1)}$)] $|\lambda_j| + |\lambda_k| \leq \lambda_c$ 
for $j\not=k$.
\item[\rm($C_J^{(1)}$)] $|\lambda_j| \leq \lambda_c/2$ 
for $j \in J$. 
\item[\rm($D_J^{(1)}$)] $|\lambda_j| + |\lambda_k| \leq \lambda_c$ 
for $j\not=k$.
\item[\rm($B_J^{(2)}$)] $|\lambda_j| \leq \lambda_c$ 
for $j\in J$. 
\item[\rm($C_J^{(2)}$)] $|\lambda_j| + |\lambda_k| \leq \lambda_c$ 
for $j\not=k$.
\item[\rm($BC_J^{(2)}$)] $|\lambda_j| \leq \lambda_c/2$ 
for $j \in J$. 
\end{description}
\end{thm}

\begin{prf} In the untwisted case $r = 1$ we have
$\hat\Delta = \{0\} \times \Delta \times \Z = X_J^{(1)}$, so that 
Proposition~\ref{prop:4.2} asserts that 
 $\lambda$ is $d$-minimal if and only if 
$|\lambda(\check\alpha)| \frac{(\alpha,\alpha)}{2}\leq \lambda_c$ for $\alpha \in \Delta.$ 

\nin $A_J^{(1)}$:  For the root system $\Delta = A_J$,  
all roots $\alpha$ satisfy $(\alpha,\alpha) = 2$, 
so that the $d$-minimality condition on $\lambda$ is 
\[ \lambda_j - \lambda_k \leq \lambda_c  \quad \mbox{ for } \quad 
j\not=k \in J.\] 
This can also be written 
as $\max\oline\lambda - \min \oline\lambda \leq  \lambda_c.$

\nin $B_J^{(1)}$: For $\Delta = B_J$, 
the roots $\eps_j$ satisfy $(\eps_j, \eps_j) = 1$ and $\check\eps_j =2 e_j$. This leads to the $d$-minimality 
conditions 
\[ |\lambda_j| \leq \lambda_c \quad \mbox{ and } \quad 
|\lambda_j \pm \lambda_k| \leq \lambda_c\]  
which is equivalent to 
$|\lambda_j| + |\lambda_k| \leq \lambda_c$ 
for $j\not=k.$

\nin $C_J^{(1)}$: For the root system $C_J$, 
the roots $2\eps_j$ satisfy $(2\eps_j)\,\check{} = e_j$ and 
$(2\eps_j, 2 \eps_j) = 4$. The $d$-minimality thus implies 
$|\lambda_j| \leq \lambda_c/2.$ 
This also implies that
$|\lambda_j\pm \lambda_k| \leq \lambda_c$ for 
$j\not=k \in J$, so that it characterizes the $d$-minimal weights. 

\nin $D_J^{(1)}$: For the root system $D_J$, 
we find the conditions 
$|\lambda_j \pm \lambda_k| \leq \lambda_c$ 
which are equivalent to 
$|\lambda_j| + |\lambda_k| \leq \lambda_c$ for $j \not=k \in J.$

\nin $B_J^{(2)}$: In this case $\Delta_0 = B_J$ and 
$\Delta_1 = \{ \pm \eps_j \: j \in J  \}$ with $\|\eps_j\| = 1$. 
In view of $\check \eps_j = 2e_j$, we obtain from the roots in $\Delta_1$ the condition 
$|\lambda_j| = \shalf |2\lambda_j| \leq \lambda_c$. 
For the roots $\eps_j \pm \eps_k \in \Delta_0$ we obtain the additional condition 
$|\lambda_j \pm \lambda_k| \leq 2 \lambda_c$ which is redundant. 

\nin $C_J^{(2)}$: In this case $\Delta_1 = D_J$ and $\Delta_0 = C_J$ 
with $\|2\eps_j\|^2 = 4$ lead to the conditions 
\[ |\lambda_j \pm \lambda_k| \leq \lambda_c \quad \mbox{ and } \quad 
|\lambda_j|  \leq \lambda_c,\] 
which is equivalent to 
$|\lambda_j| + |\lambda_k| \leq \lambda_c$ for $j\not=k \in J.$

\nin $BC_J^{(2)}$: Here $\Delta_1 = BC_J$ and $\Delta_0 = B_J$ with 
$\|2 \eps_j\|^2 = 4$ lead to the conditions 
$|\lambda_j| \leq \lambda_c/2$ for the roots $\alpha = \pm 2 \eps_j$,  
and the roots $\alpha = \pm \eps_j$ provide no additional restriction. 
For the roots $\alpha = \eps_j \pm \eps_k$ we obtain 
$|\lambda_j \pm \lambda_j|\leq \lambda_c$, which also is redundant. 
\end{prf}

\begin{rem} \mlabel{rem:4.5} (a) The preceding theorem implies that 
$d$-minimal weights $\lambda \in \hat V^*$ define bounded 
functions $\oline\lambda \: J \to \R$ and, moreover, that the boundedness of $\oline\lambda$ is 
equivalent to the existence of a $\lambda_c > 0$ such that 
$\lambda = (\lambda_c, \oline\lambda, \lambda_d) \in \hat V^*$ is $d$-minimal. 

(b) If $\lambda  \in V^*$ satisfies 
$\lambda(\check \alpha) \in \Z$ for each $\alpha \in \hat\Delta$, 
then the subset 
$\lambda + \hat\cQ \subeq V^*$, where  
$\hat\cQ = \la \hat\Delta \ra_{\rm grp}$  is the {\it root group}, 
is invariant under the Weyl group $\hat\cW$. 
Therefore $(\hat\cW\lambda)(d) 
\subeq \lambda(d) + \Z$. 
If $(\hat\cW\lambda)(d)$ is bounded from below, we thus 
obtain the existence of a $d$-minimal element in $\hat\cW\lambda$. 
\end{rem}

For general functionals which are not integral weights, 
the situation is more complicated, as 
Example~\ref{ex:4.7} below shows.

\subsection{The affine Weyl group} 
\mlabel{subsec:3.4} 

Recall the inclusion 
\[ \iota_\cW \: \cW \cong \la r_{(0,\alpha,0)} 
\: \alpha \in \Delta_0\ra_{\rm grp} \into \hat \cW\]  
of the locally finite Weyl group $\cW$ into $\hat\cW$ 
and note that it provides a section of the quotient homomorphism 
$q \: \hat \cW \to \cW$ corresponding to the passage from $\hat V$ to $V$. 
For $n \in \Z$ and $(0,\alpha,n)\in \hat\Delta$, the elements $r_{(0,\alpha,0)} r_{(0,\alpha,n)}$ 
generate the normal subgroup $\cN := \ker q$. 

To make the structure of $\cN$ more explicit, we consider for $x \in V$ the endomorphism 
$\tau_x = \tau(x)$ of $\hat V$, defined by 
\[ \tau_x(z,y,t) := 
\Big(z + \la y,x\ra + \frac{t\|x\|^2}{2}, 
y + t x, t\Big). \] 
The maps $\tau_x$ are isometries with respect to the Lorentzian form and 
an easy calculation shows that 
\[ \tau_{x_1} \tau_{x_2} = \tau_{x_1 + x_2}\quad \mbox{ for } 
\quad x_1, x_2 \in V\] 
and that 
$r_{(0,\alpha,0)} r_{(0,\alpha,n)}  = \tau_{n\check \alpha}$  for $(0,\alpha,n) \in \hat\Delta$.  This leads to 
\[ \cN =\tau(\cT) \quad \mbox { for } \quad 
\cT := \la n \check \alpha \: \alpha \in \Delta_n^\times, n \in \N\ra_{\rm grp}.\] 

\begin{prop} \mlabel{prop:4.8} 
For the untwisted root systems of type $X_J^{(1)}$, 
the group $\cT$ coincides with the group 
$\check\cR := \la \check \Delta \ra_{grp}$ 
of coroots. For the three twisted cases, it is given in 
$V$ in terms of the canonical basis elements $(e_j)_{j \in J}$ by: 
\begin{description}
\item[\rm(i)] $\cT = 2 \Z^{(J)}$ for $B_J^{(2)}$. 
\item[\rm(ii)] $\cT = \big\{ \sum_{j \in J} n_j e_j \: 
\sum_j n_j \in 2 \Z \big\}$ for $C_J^{(2)}$. 
\item[\rm(iii)] $\cT = \sum_J \Z e_j \cong \Z^{(J)}$ for $BC_J^{(2)}$. 
\end{description}
\end{prop}

\begin{prf} (i) For $B_J^{(2)}$ we derive from 
$\Delta_1 = \{ \pm \eps_j \: j \in J  \}$, 
$\Delta_0 =B_J$ and 
$\check \eps_j = 2 e_j$, 
$(\eps_j \pm \eps_k)\,\check{}\, = e_j \pm e_k$ that 
$\cT$ is the subgroup of $V$ generated by the elements 
\[ 2(e_j \pm e_k),\ j\not=k \quad \mbox{ and } \quad 2 e_j, \ j \in J.\] 

\nin (ii) For $C_J^{(2)}$ we have $\Delta_1 = D_J$ and $\Delta_0 = C_J$, 
which leads to the generators 
\[ \pm e_j \pm e_k,\ j\not=k \quad \mbox{ and } \quad 2 e_j, \ j \in J.\] 

\nin (iii) For $BC_J^{(2)}$ we obtain from $\Delta_1 = BC_J$ and 
$\Delta_0 = B_J$ the generators $\pm e_j, \pm e_j \pm e_k$ for $j \not=k$. 
\end{prf}

Note that 
\[ \hat\cW d = \cN d  
= \Big\{ \Big(\frac{\|x\|^2}{2}, x, 1\Big) \: x \in \cT\Big\}\]  
leads to 
\begin{equation}
  \label{eq:nvals}
(\hat\cW \lambda)(d) 
= \Big\{ \lambda_c \frac{\|x\|^2}{2} + \oline\lambda(x) + \lambda_d \: 
x \in \cT\Big\}.
\end{equation}
This formula shows immediately that if 
$\lambda_c > 0$ and 
\[ \|\oline\lambda\|_2^2 := \sum_{j \in J} |\lambda_j|^2 < \infty,\]  
which is in particular the case if $\supp(\oline\lambda)$ is finite, 
then $(\hat\cW\lambda)(d)$ is bounded from below. 
Since $\cT$ is not a vector space, we cannot expect 
that the condition that  $(\hat\cW\lambda)(d)$ is bounded from below  
implies that $\|\oline\lambda\|_2 < \infty$, and 
Theorem~\ref{thm:dmin} does indeed show that this is not the case. 
It only implies that  $\oline\lambda$ is bounded. 
The following example illustrates the situation further.

\begin{ex} \mlabel{ex:4.7}  
We provide an example of an element 
$\lambda \in \hat V^*$ for which 
$(\hat\cW \lambda)(d)$ is  bounded from below, but 
contains no minimum. 

We consider the root system $A_\N$ ($J = \N$) and 
$\lambda = (1,\oline\lambda,0) \in \hat V^*$ defined by 
\[ \oline\lambda \:  \N \to \R, \quad 
\lambda_{2k} = 0 \quad \mbox{ and } \quad 
\lambda_{2k-1} = 1 + \frac{1}{k^2} \quad \mbox{ for } \quad k \in \N.\] 
On $\cT = \{ x \in \Z^{(J)} \: \sum_n x_n = 0\}$ we then consider 
the function 
\[ f \: \cT \to \R, \quad 
f(x) := \shalf \|x\|^2 \lambda_c + \oline\lambda(x) 
= \shalf \|x\|^2 + \sum_{n = 1}^\infty x_n \lambda_n 
:= \shalf \|x\|^2 + \sum_{k = 1}^\infty x_{2k-1}\Big(1 + \frac{1}{k^2}\Big).\] 
We claim that $f$ is bounded from below but that it does not have a 
minimal value.

If $x_{2k-1} \leq -3$ for some $k$, then we consider the element 
$\tilde x := x + e_{2k-1} - e_{2\tilde k}$, where 
$\tilde k$ is such that $x_{2\tilde k} = 0$.  Then 
\[ f(x) - f(\tilde x) 
= \shalf(x_{2k-1}^2 - (x_{2k-1}+1)^2 - 1) - \Big(1 + \frac{1}{k^2}\Big)
= -x_{2k-1} -2 - \frac{1}{k^2} \geq 0.\] 
To show  that $f(\cT)$ is bounded from below, it therefore suffices to 
consider $f(x)$ for elements $x \in \cT$ satisfying 
$x_{2k-1} \geq -2$ for every $k \in \N$. 
This leads to 
\[ f(x) 
= \shalf \|x\|^2 + \sum_{k = 1}^\infty x_{2k-1}\Big(1 + \frac{1}{k^2}\Big)
\geq \shalf \|x\|^2 + \sum_{k = 1}^\infty x_{2k-1}
-2 \sum_{k = 1}^\infty \frac{1}{k^2},\] 
and since 
$\shalf \|x\|^2 + \sum_{k = 1}^\infty x_{2k-1} \geq 0$ for every 
$x \in \cT$ by Theorem~\ref{thm:dmin}($A_J^{(1)}$), 
we see that $f$ is bounded from below. 

If $f(x_0) = \min f(\cT)$, then 
\[ f(x_0 + x) = 
\shalf \|x+x_0\|^2 + \la \oline\lambda, x + x_0 \ra 
= \shalf \|x\|^2 + \la \oline\lambda + x_0, x \ra + f(x_0) \] 
implies that the function 
$\oline \mu := \oline\lambda + x_0$ defines a $d$-minimal functional $\mu 
= (1, \oline\mu, 0)$, so that 
\[ \sup\oline \mu - \inf\oline \mu  \leq \mu_c = 1.\] 
Since $x_0$ has finite support, 
$\sup \oline\mu> 1$, so that $\inf \uline\mu \leq 0$ leads to 
a contradiction. Therefore 
$\hat\cW\lambda$ contains no $d$-minimal element. 
\end{ex}


\begin{thebibliography}{aaaaaaa} 

\bibitem[HoG99]{HoG99} Hofmann, G., ``Invariante konvexe Mengen f\"ur lineare 
Coxetergruppen,'' Diplomarbeit, TU Darmstadt, 1999

\bibitem[Hu92]{Hu92} Humphreys, J. E., 
``Reflection Groups and Coxeter Groups,''
Cambridge studies in advanced mathematics {\bf 29}, Cambridge
University Press, 1992 

\bibitem[Ka90]{Ka90} Kac, V., ``Infinite-dimensional Lie Algebras," Cambridge University Press, $3^{rd}$ printing, 1990  

\bibitem[LN04]{LN04} Loos, O., and E. Neher, ``Locally Finite Root Systems,'' 
Memoirs of the Amer. Math. Soc., Vol. 171, {\bf 811}, 2004 

\bibitem[MP95]{MP95} Moody, R., and A. Pianzola, ``Lie Algebras with Triangular
Decompositions,'' Canad. Math. Soc. Series of Monographs and advanced texts, 
Wiley In\-ter\-science, 1995

\bibitem[Ne98]{Ne98} Neeb, K.-H., {\it Holomorphic highest weight representations
of infinite dimensional complex classical groups}, 
J.\ reine angew. Math.\ {\bf 497} (1998), 171--222  

\bibitem[Ne10]{Ne10} ---, {\it Unitary highest weight modules of 
locally affine Lie algebras}, in 
``Proceedings of 
the Workshop on Quantum Affine Algebras, Extended Affine Lie Algebras
and Applications (Banff, 2008)'', Eds. Y. Gao et al, Contemporary Math. 
{\bf  506}, Amer. Math. Soc., 2010; 227--262

\bibitem[Ne12]{Ne12} ---, {\it Semibounded representations 
of double extensions of Hilbert loop groups}, Preprint, 2012 

\bibitem[NS01]{NS01} Neeb, K.--H., and N.\ Stumme, {\it The classification of
locally finite split simple Lie algebras}, J. Algebra {\bf 553} (2001), 
25--53 

\bibitem[Vin71]{Vin71} Vinberg, E. B., {\it Discrete linear groups generated by
reflections}, Math. USSR-Izv. {\bf 5} (1971), 1083--1019 

\bibitem[YY10]{YY10} Y. Yoshii, {\it Locally extended affine root 
systems}, in this volume: ``Quantum affine algebras, extended affine Lie 
algebras and applications'', Eds. Y. Gao et al, Contemporary Math. 
{\bf  506}, Amer. Math. Soc., 2010; 285--302 


\end{thebibliography}
\end{document}